\documentclass[12pt]{article}
\usepackage{amsfonts}
\usepackage{amsmath, amssymb}
\usepackage{pifont}
\usepackage{natbib}
\usepackage{color}

\usepackage{url}

\title{Probability, Statistics and Planet Earth, \\ II:The Bochner-Godement theorem for symmetric spaces}
\author{N. H. Bingham and Tasmin L. Symons}

\begin{document}
\def\R{\mathbb{R}}
\def\C{\mathbb{C}}
\def\Z{\mathbb{Z}}
\def\N{\mathbb{N}}
\def\Q{\mathbb{Q}}
\def\D{\mathbb{D}}
\def\Sp{{\mathbb{S}}}
\def\T{\mathbb{T}}
\def\hb{\hfil \break}
\def\ni{\noindent}
\def\i{\indent}
\def\a{\alpha}
\def\b{\beta}
\def\e{\epsilon}
\def\d{\delta}
\def\De{\Delta}
\def\g{\gamma}
\def\qq{\qquad}
\def\L{\Lambda}
\def\E{\cal E}
\def\G{\Gamma}
\def\F{\cal F}
\def\K{\cal K}
\def\A{\cal A}
\def\B{\cal B}
\def\M{\cal M}
\def\P{\cal P}
\def\Om{\Omega}
\def\om{\omega}
\def\s{\sigma}
\def\th{\theta}
\def\Th{\Theta}
\def\z{\zeta}
\def\p{\phi}
\def\m{\mu}
\def\n{\nu}
\def\Si{\Sigma}
\def\q{\quad}
\def\qq{\qquad}
\def\half{\frac{1}{2}}
\def\hb{\hfil \break}
\def\half{\frac{1}{2}}
\def\pa{\partial}
\def\r{\rho}
\def\Spd{{{\Sp}}^d}

\maketitle

\ni {\bf Abstract} \\
\i The Bochner-Godement theorem extends the classical Bochner and Bochner-Schoenberg theorems from the context of Euclidean spaces and spheres to general symmetric spaces.  We show how it also includes recent results on products of symmetric spaces: the Berg-Porcu theorem (sphere cross line), and the Guella-Menegatto-Peron theorem (products of spheres), and discuss related results and applications. \\

\ni MSC: 60G60, 62H11. \\
\ni Keywords: {\sl spatio-temporal covariances; random fields; geodesic distance.}\\

\ni { \bf 1. Introduction; the Bochner-Schoenberg theorem} \\

Our motivation here is a group of three theorems, one classical and two recent, all visibly related and all crucially relevant to our motivating area, the mathematics, probability and statistics of Planet Earth (see [BinS], referred to below as I).  

We write ${\Sp}^d$ for the $d$-sphere (a $d$-dimensional manifold embedded in $(d+1)$-space ${\R}^{d+1}$).  We are interested in {\it positive definite (pd)} functions $f: {\Sp}^d \times {\Sp}^d \to \R$; we restrict attention to the {\it isotropic} case, where $f$ is of the form 
$$
f({\bf x},{\bf y}) = g(d({\bf x},{\bf y})) \qquad ({\bf x},{\bf y} \in {\Sp}^d)
$$ 
with $d$ geodesic distance on the sphere; here $g$ is called the {\it isotropic part} of $f$.  For convenience, we adopt the usual normalisation of taking the sphere to have radius 1; thus $f$ is defined on $\Sp \times \Sp$, $g$ on $[-1,1]$, and $x = d({\bf x}, {\bf y}) \in [-1,1]$.  The key role here is played by the {\it spherical functions} (the name derives from spherical harmonics: I \S 8.1, [BinMS], [Hel1,3]).  The relevant $g$ here are orthogonal polynomials (for background on which see e.g. [Sze]), the {\it Gegenbauer} (or {\it ultraspherical}) {\it polynomials}, of index $\lambda$,
$$
P_n(x), \qquad x = d({\bf x}, {\bf y}), \qquad \lambda = \half (d-1),
$$
on ${\Sp}^d \subset {\R}^{d+1}$ as above.  For the 2-sphere in 3-space (i.e. the Earth), these reduce ($\lambda = \half$) to the familiar {\it Legendre polynomials}, $P_n(x)$, and for the circle (1-sphere) in the plane to the {\it Tchebycheff polynomials}.  All these are orthogonal polynomials with orthogonality interval $[-1,1]$ [Sze]. 

The statistical relevance of positive definite functions (or by restriction, matrices) here is that they correspond to covariance functions or matrices of vectors (via Gram matrices, or Gramians); similarly for correlation functions or matrices, for vectors on the unit sphere; see e.g. [HorJ1, Ch. 7]. \\ 

\ni {\bf  Bochner-Schoenberg Theorem}.  For ${\Sp}^d$, the general isotropic positive-definite function, equivalently, the general isotropic covariance, is given (to within a scale factor) by {\it a convex combination of Gegenbauer polynomials} $P_n^{\lambda}(x)$:
$$
c \sum_{n=0}^{\infty} a_n P_n^{\lambda}(x), \quad a_n \geq 0, \quad \sum_{n=0}^{\infty} a_n = 1, \quad x = d({\bf x}, {\bf y}), \ {\bf x}, {\bf y} \in {\Sp}^d,                     \eqno(BS)
$$
with Legendre polynomials when $d = 2$. \\

\i In probabilistic language, this says that the general case is a {\it mixture} of spherical harmonics: that is, that {\it the Fourier-Gegenbauer coefficients} $a_n$ here {\it are non-negative}.

The study of positive definite functions on spheres dates back to Bochner and Schoenberg 1940-42 ([Sch1], [BocS] in 1940; Bochner [Boc2] in 1941; [Sch2] in 1942; see I for more detail).  Later extensions to probability theory were made by Gangolli [Gan], Bingham [Bin2], Faraut [Far1] and others.  The importance of this classical area has been widely recognised by the geostatistical community recently. 
     
Askey and Bingham [AskB, \S 2] extended this result to Gaussian processes on the other compact symmetric spaces of rank 1 (two-point homogenous Riemannian manifolds).  \\

\ni {\bf 2.  Symmetric spaces; the Bochner and Bochner-Godement theorems} \\

\ni {\it Bochner's theorem} \\
\i The simplest setting for a characterisation theorem for positive definite functions is the line $\R$, or $d$-space ${\R}^d$, regarded as a topological group under addition.  These were characterised by {\it Bochner's theorem} [Boc1] of 1933 as the Fourier-Stieltjes transforms of (positive, finite) measures -- in probabilistic language, as (to within a normalisation factor) {\it characteristic functions}.  The key role here is played by the {\it characters} (homomorphism) ($t \mapsto e^{itx}$ for the reals).  The theory was extended to locally compact abelian groups by Weil in 1940 [Wei], and to (not necessarily abelian) locally compact groups by Godement [God1] in 1948. \\

\ni {\it Spheres}. \\
\i The sphere ${\Sp}^d$ above is not a group but a coset space, under the identification
$$
{\Sp}^d = SO(d+1)/SO(d)
$$
with $SO(d)$ the special orthogonal group in $d$-space.  It carries an {\it involutive isometry}, where one fixes a point (say, the North Pole $O$), and maps each point $P$ into its `mirror image' $P'$ in $O$ in the geodesic linking $P$ to $O$ and produced an equal distance beyond $O$.  This gives the prototype for the compact case of a {\it symmetric space}, below.\\

\ni {\it Symmetric spaces} \\
\i A symmetric space $X$ is a metric space carrying an involutive isometry as above; these arise as coset spaces (or homogeneous spaces) $X = G/K$; the pair $(G,K)$, also spoken of as a homogeneous space, has $G$ a Lie group and $K$ a compact subgroup.  These arise from the work of Elie Cartan and others in the 1920s, in the context of Lie theory; for background and details see e.g. the books by Helgason [Hel1,2,3,4], Wolf [Wol1,2] and Faraut [Far3].  The role of the characters here is played by the {\it spherical functions}.  These were developed by Godement in 1952 [God3] by general methods, following earlier developments using Lie theory; for a modern treatment see the books cited above. \\

\ni {\it Harmonic analysis on symmetric spaces}. \\
\i The key spherical functions here are the positive definite ones.  A spherical function need not be positive definite in general.  Those that arise in the context of spheres via Gegenbauer polynomials are [Sch1]; more generally, so are those that arise in compact symmetric spaces ([Hel1, X Th. 6.5], [Hel4, IV Th. 4.2]).  The class of positive definite spherical functions is called the {\it spherical dual}, which we will denote by $\L$.  That Bochner's theorem may be extended beyond the group case above to symmetric spaces goes back to Godement [God4] in 1957 (`Plancherel case'), giving the Bochner-Godement theorem below. \\

\ni {\it Bochner-Godement theorem} \\
\i In its modern formulation, this very useful result is as follows:  \\

\ni {\bf Bochner-Godement theorem}. \\
\i  The general positive definite function $\psi$ on a symmetric space is given (to within scale) by a mixture of positive definite spherical functions ${\phi}_{\lambda}$ over the spherical dual $\L$ by a probability measure $\mu$:
$$
\psi(x) = c \int_{\L} {\phi}_{\lambda}(x) \mu(d \lambda).      \eqno(BG)
$$

For background and details, see e.g. [Wol1, Th. 9.3.4], Faraut [Far2, Th. 1.2], Faraut and Pevzner [FarP]; cf. Faraut and Harzallah [FarH, Th. 3.1], Askey and Bingham [AskB]; see also Thomas [Tho].  Taking $K = \{e\}$, the trivial subgroup, the Bochner-Godement theorem reduces to the group case of Godement, including Bochner's theorem in the Euclidean case; for the sphere, it reduces to the Bochner-Schoenberg theorem. \\

 \ni {\bf 3. Products; the Berg-Porcu and Guella-Menegatto-Peron theorems} \\
 
\ni {\it Properties of positive definite functions} \\
\i Recall some basic closure properties of the class ${\cal P} (M)$ of positive definite functions on a symmetric space $M$ (or semigroup, as in Berg et al. [BerCR]): ${\cal P}(M)$ is closed under positive scaling and convex combinations (both touched on in the above).  Further, ${\cal P}(M)$ is closed under pointwise products, which corresponds to tensor products (of representations); see e.g. [Hel4, III Prop. 12.4], [AskB, 5.2].  Similarly for product spaces (as here): by the Schur product theorem (\S 4.3 below -- see e.g. [HorJ1, 458]), if $f_1 \in {\cal P}(M_1)$, $f_2 \in {\cal P}(M_2)$, then $f_1 f_2 \in {\cal P}(M_1 \times M_2)$.  If the symmetric space is a direct product, the spherical dual is also a direct product:
$$
(G, K) = (G_1, K_1) \times (G_2, K_2),
\qquad
\L = {\L}_1 \times {\L}_2.
$$ 

\ni {\it Sphere cross line} \\
\i We turn now from sphere $\Sp^d$ to sphere cross line $\Sp^d \times \R$.  As before, we confine attention to the isotropic case for the space variable on the sphere, and similarly, we restrict to the stationary case for the time-variable on the line.  As above, we know ${\cal P}(\Sp^d)$ from the Bochner-Schoenberg theorem, and ${\cal P}(\R)$ (or ${\cal P}({\R}_+)$) from Bochner's theorem.  Taking products gives elements of ${\cal P}(\Sp^d \times \R)$, but these are {\it separable} [Gne], with space and time independent.  The question was raised by Mijatovi\'c and the authors [BinMS, \S 4.4] in 2016 of extending the results of Askey and Bingham [AskB] to sphere cross line, and answered by Berg and Porcu [BerP] in 2017. \\

\ni {\it The Berg-Porcu theorem} \\
\i We state this result below, and include its short proof, from the Bochner-Godement theorem. \\

\ni {\bf Berg-Porcu Theorem}.  The class of isotropic stationary sphere-cross-line covariances coincides with the class of mixtures of products of Gegengauer polynomials $P_n^{\lambda}(x)$ and characteristic functions ${\phi}_n(t)$ on the line:
$$
c \sum_{n=0}^{\infty} a_n P_n^{\lambda}(x) {\phi}_n(t), \qquad a_n \geq 0, \quad \sum_{n=0}^{\infty} a_n = 1. \eqno(BP)
$$

\ni {\it Proof}.  In the notation above, with the first factor the sphere ${\Sp}^d$, ${\L}_1$ is the Gegenbauer polynomials $P_n$; with the second factor the line, ${\L}_2$ is the set of characters, which can also be identified with the line:
$$
t \leftrightarrow e^{it.} = (x \mapsto e^{itx}).
$$
We use disintegration in $(BG)$ (Fubini's theorem extended beyond product measures: see e.g. Kallenberg [Kal, Th. 6.4]), integrating the probability measure $\mu$ on $\Lambda$ first over the $x$-variable above for fixed $n$.  This gives a probability measure, ${\mu}_n$ say; integrating the character $e^{itx}$ over ${\mu}_n(dx)$ gives its characteristic function ${\phi}_n(t)$, the second factor in $(BP)$; the remaining integration is a summation over $n$, giving the first factor $P_n^{\lambda}$ in $(BP)$.  Equivalently, one can take $\lambda = ({\lambda}_1, {\lambda}_2)$ in $(BG)$ as a random variable with law $\mu$, condition on its second coordinate, and use the Conditional Mean Formula [Wil2, 390] (a special case of the tower property [Wil1, 9.7i] or chain rule [Kal, 105]).    \hfill $\square$ \\

\ni {\it The Guella-Menegatto-Peron theorem} \\
\i A second important case is that of a product of spheres.  This is the context of the Guella-Menegatto-Peron theorem below [GueMP].  The proof is immediate from the Bochner-Godement theorem. \\

\ni {\bf Guella-Menegatto-Peron theorem}.  The general isotropic positive definite function on ${\Sp}^{d_1} \times {\Sp}^{d_2}$ is of the form
$$
c \sum_{m,n 
= 0}^{\infty} a_{mn} P_m^{{\lambda}_1}(x_1) P_n^{{\lambda}_2}(x_2), 
\quad a_{mn}\geq 0, \quad \sum_{m,n = 0}^{\infty} a_{mn} = 1,
$$
$$ 
x_i = d({\bf x}_i, {\bf y}_i), 
\quad {\bf x}_i, {\bf y}_i \in {\Sp}^{d_i}, \quad{\lambda}_i = \half (d_i - 1).
$$

The function here factorises into a product of functions of $x_1$ and $x_2$ when the matrix $(a_{mn})$ is of rank one (that is, of the form $a_{mn} = b_m c_n$); in probabilistic language, this corresponds to independence. 

Multiple products over any number of spheres and lines follow from the Bochner-Godement theorem in the same way. 

Following [BerP] and [GueMP], three of the authors returned to the area [BerPP], using similar methods. \\

\ni {\bf 4.  Remarks} \\

\ni {\it 4.1.  The Gelfand-Naimark-Segal (GNS) construction}. \\
\i For $G$ a topological group, $K$ a closed subgroup, $\pi$ a unitary representation of $G$ in a Hilbert space $H$, a vector $u \in H$ is {\it cyclic} if $\{ \pi (g) u: g \in G \}$ generates a dense subspace of $H$.  The Gelfand-Naimark-Segal (GNS) construction (Gelfand and Naimark [GelN] in 1943, Segal [Seg] in 1947; [God2]) says that for $\phi \in {\cal P}(K\backslash G/K)$, there exists a unitary representation$({\pi}_{\phi}, H_{\phi})$ of $G$ and a $K$-invariant cyclic unit vector $u$ with
$$
\phi(g) = (u, \pi(g) u) \qquad (g \in G),   \eqno(GNS)
$$
and $({\pi}_{\phi}, H_{\phi}, u)$ is unique up to isomorphism.  For background, see e.g. Dixmier [Dix, \S 13.4.4].  \\
\i The GNS construction permeates the modern treatment of spherical functions on symmetric spaces; see e.g. [Hel1, X.4], [Hel3, IV.1, IV.3], [Hel4, III, Th. 12.1].  It is easily seen that any $\phi$ as in $(GNS)$ is positive definite (see e.g. [Hel1, 412-3], [Hel3, IV Th. 4.2]; these are the spherical functions that arise in symmetric spaces of compact type [Hel3, IV Th. 4.2] such as spheres and products of spheres. \\   

\ni {\it 4.2.  Extreme points and irreducible representations}. \\
\i In $(GNS)$ above, $\phi$ is an extreme point in the convex set ${\cal P}(K \backslash G / K)$ if and only if the unitary representation $(\pi_{\phi}, H_{\phi})$ is irreducible (see e.g. [Far2, Prop. 1.4], Simon [Sim, Ex. 8.21], [Hel4, Th. 12.6]).  This shows very clearly the role of convexity, the Krein-Milman theorem and positivity here; in particular, one sees that $(BG)$ gives a Choquet representation. \\
\i The GNS construction has its origins in quantum mechanics; the $\phi$ here are called {\it states}; the extreme $\phi$ are called {\it pure states}.  This physical interpretation is taken further by Davies [Dav]. \\ 

\ni {\it 4.3.  The Schur product theorem}. \\
\i The Schur product theorem (I. J. Schur in 1911, writing as J. Schur) states that the product of positive definite functions is positive definite.  This fundamental result underlies the treatment of product spaces above.  For details, see e.g. [Schu], [HorJ1, \S 7.5], [HorJ2, \S 5.2], [BerCR], [Hel4, Prop. 12.4]. \\

\ni {\it 4.4.  Gelfand pairs}. \\
\i For $G$ a locally compact group (Lie group will suffice for us here) with compact subgroup $K$, $(G,K)$ is called a {\it Gelfand pair} if the convolution algebra of $K$-biinvariant compactly supported continuous measures on $G$ is commutative; equivalently, if for any locally convex irreducible representation $\pi$ of $G$, its space of $K$-invariant vectors is at most one-dimensional.  For symmetric spaces $X = (G,K)$, $(G,K)$ is a Gelfand pair; see e.g. [Hel1, X, Th. 2.9, 4.1], [Hel3, IV, Th. 3.1]; as this is where our interest lies here, we have chosen to use the language of symmetric spaces rather than that of Gelfand pairs.  Nevertheless, Gelfand pairs are very interesting objects in their own right; see e.g. the monograph by van Dijk [vDij].  They occur in tree models; also, finite Gelfand pairs arise in many areas in combinatorics. \\  

\ni {\it 4.5.  Special functions and group representations}. \\
\i Special functions (Askey, a great expert, insists that they should instead be called useful functions) arise in many places relevant here.  They occur most obviously in the Gegenbauer polynomials in the Bochner-Schoenberg theorem of \S 1.  These functions generally arose in applications in the 19th century, before the development of group representations.  However, this theory provides a natural and powerful framework in which to study many of the classical special functions; for monograph treatments, see e.g. Vilenkin [Vil], Vilenkin and Klimyk [VilK], Andrews, Askey and Roy [AndAR]. \\

\ni {\it 4.6.  Spaces of constant curvature}. \\
\i The theory of symmetric spaces is simplest for those of {\it rank one} [Hel1,2,3,4].  These are simple enough to be {\it classified}.  They comprise, apart from certain exceptional cases, the {\it spaces of constant curvature}, $\kappa$.  There are three cases: \\
(i) $\kappa > 0$.  These are the {\it spheres}, in any dimension $d$; the curvature $\kappa$ depends on the radius; with our normalisation of radius 1, we obtain the ${\Sp}^d$. \\
(ii) $\kappa = 0$.  These are the {\it Euclidean spaces}, again in any dimension $d$. \\
(iii) $\kappa < 0$.  These are the {\it hyperbolic spaces}.\\
For background and details, see e.g. Wolf [Wol2], [Hel1, X], [Hel3, IV, V]. \\

\ni {\it 4.7.  Hypergroups}. \\
\i The spaces of constant curvature provide concrete and accessible examples of {\it hypergroups}.  These are spaces endowed with a convolution-like operation satisfying certain properties; the concept emerged from work of Dunkl, Jewett and Spector; for a detailed treatment and a wealth of examples, see the standard work on the subject, Bloom and Heyer [BloH]. \\
\i It turns out that the spaces of constant curvature provide probabilistically interesting examples of hypergroups precedng the hypergroup concept.  The spheres give the {\it Bingham} (or {\it Bingham-Gegenbauer}) {\it hypergroup} ([Bin1], [BloH]); the Euclidean case gives the {\it Kingman} (or {\it Kingman-Bessel}) {\it hypergroup} ([Kin], [BloH]); the hyperbolic spaces give the {\it hyperbolic hypergroup}.\\
 
\ni {\it 4.8.  Connections with number theory}. \\
\i Particularly in the hyperbolic case (\S 4.6), there are many connections with number theory -- for example, with the Selberg trace formula and its relatives (this explains the title of [God4]).  For background and references here, see e.g. the books of Terras [Ter1,2]. \\     

\ni {\it 4.9.  Multivariate applications}. 

The results discussed above extend to the multivariate setting. In this case, one needs to consider {\sl cross-covariance} functions. Recent work [AlePFM] extends the traditional Euclidean framework to the $\Sp^d \times \R$ setting discussed here, offering some parametric families of matrix-valued covariances on sphere cross line.\\

\ni {\it 4.10.  Anisotropy}. 

Throughout, we have assumed the spatial component of the process is {\sl isotropic} -- the covariance between two points $x$ and $y$ depends solely on the distance between them. This is a convenient assumption -- indeed, it facilitates all the theory exposited above -- and is pervasive in spatial statistics: see [Ma] for a discussion of isotropy's impact on covariance modelling over $\R^d$.\\

\ni {\bf Acknowledgement}. \\

\i As in I, the second author thanks the EPSRC Mathematics of Planet Earth CDT, based jointly in the Mathematics Department at Imperial College, London and the School of Mathematical and Physical Sciences at the University of Reading, for financial support during her doctoral studies.  Both authors thank the same sources for academic support. \\  

\ni{\bf References} \\

\ni [AlePFM] A. Alegria, E. Porcu, R. Furrer and J. Mateu. Covariance functions for multivariate Gaussian fields evolving temporally over Planet Earth. {\url{arXiv: 1701.06010v1}}.\\
\ni [AndAR] G. E. Andrews, R. Askey and R. Roy, {\sl Special functions}.  Encycl. Math. Appl. {\bf 71}, Cambridge University Press, 1999. \\
\ni [AskB] R. A. Askey and N. H. Bingham, Gaussian processes on compact symmetric spaces.
{\sl Z. Wahrschein.  verw. Geb.} {\bf 37} (1976), 127-143. \\
\ni [BerCR] C. Berg, J. P. R. Christensen and P. Ressel, {\sl Harmonic analysis on semigroups: Theory of positive definite and related functions}.  Springer, 1984. \\
\ni [BerP] C. Berg and E. Porcu, From Schoenberg coefficients to Schoenberg functions. {\sl Constructive Approximation} {\bf 45} (2017), 217 -- 241.\\
\ni [BerPP] C. Berg, A. P. Peron and E. Porcu, Orthogonal expansions related to compact Gelfand pairs.  arXiv:1612.03718v1.  \\
\ni [Bin1] N. H. Bingham, Random walk on spheres.  {\sl Z. Wahrschein. Verw. Geb.} {\bf 22} (1972), 169-192. \\
\ni [Bin2]  N. H. Bingham, Positive definite functions on spheres.  {\sl Proc. Cambridge Phil. Soc.} {\bf 73} (1973), 145-156. \\
\ni [BinMS] N. H. Bingham, Aleksandar Mijatovi\'c and Tasmin L. Symons. Brownian manifolds, negative type and geo-temporal covariances. {\sl Communications in Stochastic Analysis (H. Heyer Festschrift)}, {\bf 10}(4) (2016), 421 -- 432 (arXiv:1612.06431). \\
\ni [BinS] N. H. Bingham and Tasmin L. Symons, Probability, statistics and Planet Earth, I: Geotemporal covariances, arXiv:1706.02972v2. \\
\ni [BloH] W. R. Bloom and H. Heyer, {\sl Harmonic analysis of probability measures on hypergroups}.  De Gruyter Studies in Math. {\bf 20}, Walter de Gruyter, Berlin, 1994. \\
\ni [Boc1] S. Bochner, Monotone Funktionen, Stieltjessche Integrale und harmonische Analyse.  {\sl Math. Annalen} {\bf 108} (1933), 378-410. \\
\ni [Boc2] S. Bochner, Hilbert distances and positive definite functions.  {\sl Ann. Math.} {\bf 42} (1941), 647-656.  \\
\ni [BocS] S. Bochner and I. J. Schoenberg, On positive definite functions on compact spaces.  {\sl Bull. Amer. Math. Soc.} {\bf 46} (1940), 881. \\
\ni [Dav] E. B. Davies, Hilbert-space representations of Lie algebras.  {\sl Comm. Math. Phys.} {\bf 23} (1971), 159-168. \\
\ni [vDij] G. van Dijk, {\sl Introductin to harmonic analysis and generalized Gelfand pairs}.  De Gruyter Studies in Math. {\bf 36}, Walter de Gruyter, 2009. \\
\ni [Dix] J. Dixmier, $C^{\ast}$-{\sl algebras}.  North-Holland, 1982. \\
\ni [Far1] J. Faraut, Fonction brownienne sur une vari\'et\'e Riemannienne.  {\sl S\'em. Prob. VII}, 61-76, Lecture Notes in Math.  {\bf 321} (1973). \\
\ni [Far2] J. Faraut, {\sl Infinite-dimensional spherical analysis}. COE Lecture Note {\bf 10}, Kyushu University, 2008 (Notes by Sho Matsumoto). \\
\ni [Far3] J. Faraut, {\sl Analysis on Lie groups}.  Cambridge Studies in Advanced Math. {\bf 110}, Cambridge University Press, 2008. \\
\ni [FarH] J. Faraut and K. Harzallah, Distances hilbertiennes invariantes sur un espace homog\`ene.  {\sl Ann. Inst. Fourier} {\bf 24}.3 (1974), 171-283. \\
\ni [FarP] J. Faraut and M. Pevsner, Berezin kernels and analysis in Makarevich spaces.  {\sl Indag. Math.} {\bf 16} (2005), 461-486. \\
\ni [Gan] R. Gangolli, Positive definite kernels on homogeneous spaces and certain stochastic processes related to L\'evy's Brownian motion of several parameters.  {\sl Ann. Inst. H. Poincar\'e B (NS)} {\bf 3} (1967), 121-226. \\
\ni [GelN] I. M. Gelfand and M. A. Naimark, On the embedding of normed rings into the ring of operators in Hilbert space (Russian).  {\sl Mat. Sbornik} NS {\bf 12} (1943), 197-213.\\  
\ni [Gne] T. Gneiting, Non-separable, stationary covariance functions for space-time data.  {\sl J. Amer. Stat. Assoc.} {\bf 87} (2002), 590-600. \\
\ni [God1] R. Godement, Les fonctions de type positif et la th\'eorie des groupes.  {\sl Trans. Amer. Math. Soc.} {\bf 63} (1948), 1-84. \\
\ni [God2] R. Godement, Sur la th\'eorie des repr\'esentations unitaires.  {\sl Ann. Math.} {\bf 53} (1951), 65-120. \\
\ni [God3] R. Godement, A theory of spherical functions I.  {\sl Trans. Amer. Math. Soc.} {\bf 73} (1952), 496-556. \\
\ni [God4] R. Godement, Introduction aux travaux de A. Selberg.  {\sl S\'eminaire Bourbaki} {\bf 4} Exp. 144, 95-110, Soc. Math. France, 1957. \\
\ni [GueMP] J. C. Guella, V. A. Menegatto and A. P. Peron, An extension of a theorem of Schoenberg to products of spheres.  {\sl Banach J. Math. Analysis} {\bf 10} (2016), 671-685. \\
\ni [Hel1] S. Helgason, {\sl Differential geometry and symmetric spaces}.  Academic Press, 
1962. \\
\ni [Hel2] S. Helgason, {\sl Differential geometry, Lie groups and symmetric spaces}.  Academic Press, 1978. \\
\ni [Hel3] S. Helgason, {\sl Groups and geometric analysis: Integral geometry, invariant differential operators, and spherical functions}.  Academic Press, 1978 (2nd ed., Amer. Math. Soc., 2001). \\
\ni [Hel4] S. Helgason, {\sl Geometric analysis on symmetric spaces}.  Math. Surveys Monog. {\bf 39}, Amer. Math. Soc., 1994 (2nd ed. 2008). \\
\ni [HorJ1] R. A. Horn and C. R. Johnson, {\sl Matrix analysis}, 2nd ed., Cambridge University Press, 2013. \\
\ni [HorJ2] R. A. Horn and C. R. Johnson, {\it Topics in matrix analysis}, Cambridge University Press, 1991. \\
\ni [Kal] O. Kallenberg, {\sl Foundations of modern probability theory}, 2nd ed., Springer, 2002 (1st ed. 1997). \\
\ni [Kin] J. F. C. Kingman, Random walks with spherical symmetry.  {\sl Acta Math.} {\bf 109} (1963), 11--53. \\
\ni [Ma] C. Ma. Why is isotropy so prevalent in spatial statistics? {\sl Proc. Amer. Math. Soc} {\bf 135}(3) (2007), 865 -- 871.\\
\ni [Sch1] I. J. Schoenberg, On positive definite functions on spheres. {\sl Bull. Amer. Math. Soc} {\bf 46} (1940), 888. \\  
\ni [Sch2] I. J. Schoenberg, Positive definite functions on spheres.  {\sl Duke Math. J. } {\bf 9} (1942), 96-108 ({\sl Selected Papers 1}, 172-185). \\
\ni [Schu] I. J. Schur, Bemerkungen zur Theorie der Bilinearformen mit unendlich vielen Ver\"anderlichen.  {\sl J. reine angew. Math.} {\bf 140} (1911), 1-100. \\
\ni [Seg] I. E. Segal, Irreducible representations of operator algebras.  {\sl Bull. Amer. Math. Soc.} {\bf 53} (1947), 73-88.  \\
\ni [Sim] B. Simon, {\sl Convexity}.  Cambridge Tracts in Math. {\bf 187}, Cambridge University Press, 2011. \\
\ni [Sze] G. Szeg\"o, {\sl Orthogonal polynomials}.  Amer. Math. Soc. Colloq. Publ. {\bf XXIII}, Amer. Math. Soc., 1939. \\
\ni [Ter1] A. Terras, {\sl Harmonic analysis on symmetric spaces -- Euclidean space, the sphere, and the Poincar\'e upper half-plane}.  Springer, 1985 (2nd ed. 2013). \\
\ni [Ter2] A. Terras, {\sl Harmonic analysis on symmetric spaces and applications}, Vol. 2, Springer, 1988. \\ 
\ni [Tho] E. G. F. Thomas, The theorem of Bochner-Schwartz-Godement for generalized Gelpairs.  {\sl Functional analysis: surveys and recent results III (Paderborn, 1983)}, 291-304.  Notas Mat. {\bf 94}, North-Holland, 1984. \\  
\ni [Vil] N. J. Vilenkin, {\sl Special functions and the theory of group representations}.  Transl. Math. Monog. {\bf 22}, Amer. Math. Soc. 1968 (Russian original, Nauka, Moscow, 1965).\\
\ni [VilK] N. J. Vilenkin and A. U. Klimyk, {\sl Representations of Lie groups and special functions}.  Kluwer, Dordrecht, 1995. \\
\ni [Wei] A. Weil, {\sl L'int\'egration dans les groupes topologiques}.  Actualit\'es Sci. Ind. 869, Hermann, Paris, 1940. \\
\ni [Wil1] D. Williams, {\sl Probability with martingales}.  Cambridge University Press, 1991. \\
\ni [Wil2] D. Williams, {\sl Weighing the odds: A course in probability and statistics}.  Cambridge University Press, 2001. \\
\ni [Wol1] Wolf, J. A.: {\sl Spaces of constant curvature}.  Amer. Math. Soc., 1967 (6th ed. 2011). \\
\ni [Wol2] Wolf, J. A.: {\sl Harmonic analysis on commutative spaces}.  Amer. Math. Soc., 2007. \\

\ni N. H. Bingham, Mathematics Department, Imperial College, London SW7 2AZ; n.bingham@ic.ac.uk\\

\ni Tasmin L. Symons, Mathematics Department, Imperial College, London SW7 2AZ; tls111@ic.ac.uk

\end{document}